\documentclass[13pt]{article}
\usepackage[cp1251]{inputenc}
\usepackage[english]{babel}
\usepackage{amssymb,amsmath,amsthm,indentfirst}
\usepackage{multicol}
\usepackage{graphicx}

\begin{document}

\begin{center}
\LARGE{Exact Controllability of the Distributed System Governed by the Wave Equation with Memory}
\end{center}

\vskip 0.5cm
\begin{center}
{\large Igor Romanov}
\footnote{National Research University Higher School of Economics,\\
20 Myasnitskaya Ulitsa, Moscow 101000, Russia}
{\footnote{
E-mail: \it{ivromm1@gmail.com}}}
{\large Alexey Shamaev}
\footnote{Institute for Problems in Mechanics RAS,\\
101 Prosp. Vernadskogo, Block 1, Moscow 119526, Russia}
\footnote{Lomonosov Moscow State University,\\
GSP-1, Leninskie Gory, Moscow 119991, Russia}

\end{center}




\begin{abstract}
We will consider the exact controllability of the distributed system governed by the wave equation with memory. It will be proved that this mechanical system can be driven to rest in finite time, the absolute value of the distributed control function being bounded. In this case the memory kernel is a linear combination of exponentials.
\end{abstract}

\par Keyword:
Controllability to rest, equation with memory, distributed control
\par MSC 2010:
93C20, 35L53

\section{Introduction}

In this article we will consider the problem of the exact controllability of a system governed by integro-differential equation
\begin{equation}
\label{1}
\theta_{tt}(t,x)-K(0)\Delta\theta(t,x)-\int\limits^t_0K^{\prime}(t-s)\Delta\theta(s,x)ds=u(t,x),\ \ x\in\Omega,\ \ t>0.
\end{equation}
\vspace{0.25mm}
\begin{equation}
\label{2}
\theta|_{t=0}=\varphi_0(x),\quad \theta_t|_{t=0}=\varphi_1(x),
\end{equation}
\vspace{0.25mm}
\begin{equation}
\label{3}
\theta|_{\partial\Omega}=0,
\end{equation}
here
$$
K(t)=\sum\limits_{j=1}^{N}\frac{c_j}{\gamma_j}e^{-\gamma_jt},
$$
where $c_j$, $\gamma_j$ are given positive constant numbers, $u(t,x)$ is a control supported (in $x$) on a bounded domain $\Omega$
and $|u(t,x)|\leqslant M$, $M>0$ is a given constant number. The goal of the control is to drive this mechanical system to rest in finite time.
We say that the system is \emph{controllable to rest} when for all initial conditions $\varphi_0$, $\varphi_1$
it is possible to find a control $u(t,x)$ and a time $T>0$ such that $u(t,x)$ is equal to zero for any $t>T$ and the corresponding solution $\theta(t,x,u)$ of problem (\ref{1})--(\ref{3}) equals zero for any $t>T$ too.

Similar problems for  membranes and plates were studied earlier in \cite{Chernousko}. It was proved that vibrations of these mechanical
systems could be driven to rest by applying bounded (in absolute value) and volume-distributed control functions. The existence of a
bounded (in absolute value) boundary control that drives a string to rest was proved in \cite{Butk}. In this case the so-called moment problem
was effectively applied. An overview of the results concerning the boundary controllability of distributed systems can be found in \cite{Lions}.
Problems of controllability of systems similar to (\ref{1}) were considered in \cite{Ivanov}. A condition was formulated under which a solution
to the heat equation with memory could not be driven to rest in a finite time. This condition is that there are roots of some analytic function of complex
variable in the domain of holomorphism.

Problems similar to (\ref{1})--(\ref{3}) for integrodifferential equations were studied earlier in many articles.
Equation (\ref{1}) was firstly derived in \cite{Gurtin}.
The solvability and asymptotic behaviour for an abstract equation of this type  were investigated,
for example in \cite{Dafermos} and \cite{Desh}. In \cite{Rivera}, it was proved that the energy for
some dissipative system decays polynomially when the memory kernel decays exponentially.
Problems of the solvability of system (\ref{1})--(\ref{3}) were considered in \cite{Vlasov} and \cite{Vlasov_1}.
It was proved that the solution belongs to some Sobolev space on the semi-axis (in $t$)
when the kernel $K(t)$ is the series of exponentials, each exponential function tending to
zero when $t\rightarrow+\infty$. The explicit formulae for the solution of (\ref{1})--(\ref{3})
were obtained in \cite{Rautian}. In this case kernel $K(t)$ is also the series of decreasing exponentials.
It follows from these formulae that solutions tend to zero when $t\rightarrow+\infty$.
In all these articles kernels in integral summands of the equation are suggested to be none-increasing.
The problem of controllability to rest for one dimensional string equation is considered in \cite{Biccaria}. In this case, the kernel in the integral term of the equation is identical to 1 and control
is focused on a compact (part of a string) that moves at a constant speed.

In this article we also consider the so-called "null controllability". It means that
for every initial conditions $\varphi_0$ and $\varphi_1$ there is a control $u(t,x)$
and a time $T>0$ such that the corresponding solution $\theta(t,x,u)$ and its first derivative, with respect to $t$
equals zero for $t=T$. Null controllability and controllability to rest are not the same for systems with memory. In many cases,
controllability to rest is impossible.
Let us consider for example the one-dimensional case ($\Omega$ is an interval $(0,\pi)$). We prove now that the system (\ref{1}) is uncontrollable to rest if $u(t,x)\in C([0,\infty), L_2(0,\pi))$ is supported (in $x$), as well as in \cite{Ivanov}, on an interval $[a,b]$ which is properly contained in $[0,\pi]$.
It means that $u(t,x)\equiv 0$ outside $[a,b]$.
It is clear that the equation (\ref{1}) can be written in the following form
$$
\frac{\partial}{\partial t}\left(\theta_{t}(t,x)-\int\limits^t_0K(t-s)\theta_{xx}(s,x)ds-\int\limits^t_0u(s,x)ds\right)=0.
$$

Obviously function $\theta(t,x)$ is a solution of the equation (\ref{1}) if and only if this function is a solution of the following equation
\begin{equation}
\label{4}
\theta_{t}(t,x)-\int\limits^t_0K(t-s)\theta_{xx}(s,x)ds-\int\limits^t_0u(s,x)ds=f(x),
\end{equation}
where $f(x)$ is an arbitrary function. Let $t$ in (\ref{4}) be equal to zero then we obtain
$$
f(x)=\varphi_1(x).
$$

Let $\varphi_1(x)\equiv 0$. We introduce
$$
P(t,x)=\int\limits_0^t u(s,x)ds.
$$
Thus the problem (\ref{1})--(\ref{3}) reduces to the form
\begin{equation}
\label{5}
\theta_{t}(t,x)-\int\limits^t_0K(t-s)\theta_{xx}(s,x)ds=P(t,x),\quad x\in(0;\pi),\quad t>0.
\end{equation}
\begin{equation}
\label{6}
\theta|_{t=0}=\varphi_0(x),
\end{equation}
\begin{equation}
\label{7}
\theta|_{x=0}=0,\quad \theta|_{x=\pi}=0.
\end{equation}

Note that $P(t,x)\equiv 0$ outside $[a,b]$ and $P(t,x)$ can be considered as a new control function. It is the problem considered in \cite{Ivanov} (in a more general case).
If $K(t)$ is a linear combination of two exponentials then the system (\ref{5})--(\ref{7}) is uncontrollable to rest. It means that
there is an initial condition $\varphi_0$ such that, for any control $P(t,x)$, where $P(t,x)$ belongs to the corresponding space, the solution
of (\ref{5})--(\ref{7}) can not be driven to rest.
Using arguments similar to the above it can be proved that the system (\ref{1})--(\ref{3}) is uncontrollable to rest if $K(t)$ is a linear
combination of $N$ exponentials, where $N\geqslant 2$.
In \cite{JOTA} this result was generalized to the multidimensional case. More accurately,
controllability to rest is impossible if $N\geqslant 2$ and control function is applied only to the sub-domain.
But, as it will be shown in this article, the system (\ref{1})--(\ref{3}) is controllable to rest
when the control is distributed on the whole domain $\Omega$.

\section{Preliminaries}

Let $A:=-\Delta$ be an operator acting on a space $D(A):=H^2(\Omega)\cap H^1_0(\Omega)$, $\Omega\subset R^s$
$(s\in\mathbf{N})$ is a bounded domain with a smooth boundary.
Let also $\{\psi_n(x)\}_{n=1}^{+\infty}$ be a corresponding orthonormal system of eigenfunctions
and $\{\alpha^2_n\}_{n=1}^{+\infty}$ are corresponding eigenvalues such as
$\Delta\psi_n(x)+\alpha^2_n\psi_n(x)=0$.

We denote $W^2_{2,\gamma}(R_+,A)$ the linear space of functions $f:R_+=(0,+\infty)\rightarrow D(A)$ equipped with the norm
$$
\|\theta\|_{W^2_{2,\gamma}(R_+,A)}=\left(\int\limits_0^{+\infty}e^{-2\gamma t}\left(\left\|\theta^{(2)}(t)\right\|^2_{L_2(\Omega)}+\left\|A\theta(t)\right\|^2_{L_2(\Omega)}
\right)dt\right)^{\frac{1}{2}},\quad \gamma\geqslant 0.
$$

\textbf{Definition 1.}
A function $\theta(t,x)$ is called a strong solution of the problem (\ref{1})--(\ref{3}) if for some
$\gamma\geqslant 0$
this function belongs to the space $W^2_{2,\gamma}(R_+,A)$, satisfies the equation (\ref{1})
nearly everywhere (in $t$)
on the positive semiaxis $R_+$ and satisfies the initial conditions (\ref{2}).

Let us denote a function of complex variable $\lambda$
$$
l_n(\lambda):=\lambda^2+\alpha_n^2\lambda\hat{K}(\lambda),
$$
where
$$
\hat{K}(\lambda)=\sum\limits_{k=1}^{N}\frac{c_k}{\gamma_k(\lambda+\gamma_k)}.
$$

Now we formulate two theorems (see \cite{Rautian}) which are dedicated to the correct solvability of the initial
boundary value problem (\ref{1})--(\ref{3}).

\textbf{Theorem 1.}
Let $u(t,x)\equiv 0$ when $t\in R_+$, the function $\theta(t,x)\in W^2_{2,\gamma}(R_+,A)$, $\gamma>0$, is a strong solution of
the problem (\ref{1})--(\ref{3}), then for any $t\in R_+$ the following representation is obtained:
$$
\theta(t,x)=\frac{1}{\sqrt{2\pi}}\sum_{n=1}^{\infty}
\frac{(\varphi_{1n}+\lambda^+_n\varphi_{0n})e^{\lambda^+_nt}\psi_n(x)}
{l^{(1)}_n(\lambda^+_n)}+
\frac{1}{\sqrt{2\pi}}\sum_{n=1}^{\infty}\frac{(\varphi_{1n}+\lambda^-_n\varphi_{0n})e^{\lambda^-_nt}\psi_n(x)}
{l^{(1)}_n(\lambda^-_n)}+
$$
\begin{equation}
\label{8}
+\frac{1}{\sqrt{2\pi}}\sum_{n=1}^{\infty}\left(\sum\limits_{k=0}^{N-1}\frac{(\varphi_{1n}-q_{k,n}
\varphi_{0n})e^{-q_{k,n}t}}{l^{(1)}_n(-q_{k,n})}\right)\psi_n(x),
\end{equation}
where $-q_{k,n}$ are real zeros of the function $l_n(\lambda)$ ($q_{0,n}=0$, $q_{k,n}>0$, $k=1,...,N-1$),
$\lambda_n^{\pm}$ is a pair of complex conjugate nulls, $l^{(1)}_n$ is the first derivative of $l_n$
and the series (\ref{8}) converges in the norm of the space $L_2(\Omega)$.

\textbf{Theorem 2.}
Let $u(t,x)\in C([0,T],L_2(\Omega))$ for any $T>0$, $\theta(t,x)\in W^2_{2,\gamma}(R_+,A)$
is a strong solution of the problem
(\ref{1})---(\ref{3}) for some $\gamma>0$, $\varphi_0=\varphi_1=0$. Then for any $t\in R_+$ the following representation is obtained:
$$
\theta(t,x)=\frac{1}{\sqrt{2\pi}}\sum_{n=1}^{\infty}\omega_n(t,\lambda^+_n)\psi_n(x)
+\frac{1}{\sqrt{2\pi}}\sum_{n=1}^{\infty}\omega_n(t,\lambda^-_n)\psi_n(x)+
$$
\begin{equation}
\label{9}
+\frac{1}{\sqrt{2\pi}}\sum_{n=1}^{\infty}\left(\sum\limits_{k=0}^{N-1}\omega_n(t,-q_{k,n})\right)\psi_n(x),
\end{equation}
where
$$
\omega_n(t,\lambda)=\frac{\int\limits_0^tu_n(s)e^{\lambda(t-s)}ds}{l^{(1)}_n(\lambda)},
$$
$u_n(t)$ is a Fourier coefficient of $u(t,x)$ and the series (\ref{9}) converges in the norm of the space $L_2(\Omega)$.

The following lemma should be stated.

\textbf{Lemma 1.}
For any natural number $n$ the equality holds
$$
\frac{1}{l_n^{(1)}(\lambda^+_n)}+\frac{1}{l_n^{(1)}(\lambda^{-}_n)}+
\sum\limits_{k=0}^{N-1}\frac{1}{l_n^{(1)}(-q_{k,n})}=0.
$$

\begin{proof}
We shall deal with the solution of the problem (\ref{1})--(\ref{3}) in the case of $\varphi_0=\varphi_1=0$.
According to the theorem 2 this solution has the form of (\ref{9}), $u(t,x)$ being arbitrary and satisfying theorem conditions.
Taking partial derivative of $\theta(t,x)$ with respect to $t$ we obtain
$$
\frac{\partial\theta(t,x)}{\partial t}=\frac{1}{\sqrt{2\pi}}\sum_{n=1}^{\infty}\left(\frac{1}{l_n^{(1)}(\lambda^+_n)}+\frac{1}{l_n^{(1)}(\lambda^{-}_n)}
+\sum\limits_{k=0}^{N-1}\frac{1}{l_n^{(1)}(-q_{k,n})}\right)u_n(t)\psi_n(x)+
$$
$$
+\frac{1}{\sqrt{2\pi}}\sum_{n=1}^{\infty}\lambda^+_n\omega_n(t,\lambda^+_n)\psi_n(x)+
\frac{1}{\sqrt{2\pi}}\sum_{n=1}^{\infty}\lambda^-_n\omega_n(t,\lambda^-_n)\psi_n(x)+
$$
\begin{equation}
\label{10}
+\frac{1}{\sqrt{2\pi}}\sum_{n=1}^{\infty}\left(\sum\limits_{k=1}^{N-1}(-q_{k,n})\omega_n(t,-q_{k,n})\right)\psi_n(x).
\end{equation}

Since $\theta_t(t,x)|_{t=0}=0$ then for any natural number $n$ from (\ref{10}) arises
\begin{equation}
\label{11}
\left(\frac{1}{l_n^{(1)}(\lambda^+_n)}+\frac{1}{l_n^{(1)}(\lambda^{-}_n)}+
\sum\limits_{k=0}^{N-1}\frac{1}{l_n^{(1)}(-q_{k,n})}\right)u_n(0)=0.
\end{equation}
 By virtue of the fact that $u(t,x)$ is arbitrary, it is chosen in such a way  that all its Fourier coefficients $u_n(t)$ with respect to $t=0$ are nonzero. Thus, dividing by $u_n(0)$ in the equation (\ref{7}), we obtain the required statement. Lemma is proved.
\end{proof}

Let us consider the space $l_{\beta}$ of sequences $\{c_n\}_{n=1}^{+\infty}$ such that the series
$$
\sum\limits_{n=1}^{+\infty}|c_n|^2\alpha^{2\beta}_n
$$
converges. Then we define the space
$$
D(A^{\frac{\beta}{2}})=\left\{f(x)=\sum\limits_{n=1}^{+\infty}f_n\psi_n(x):
\{f_n\}_{n=1}^{+\infty}\in l_{\beta}\right\}.
$$

\section{The main results}

The section is devoted to the proof of the main theorem which states that the solution and its first derivative with respect to time can be driven to the null state during finite time. Then it shows that the control function which is constructed in the proof drives the system to rest actually.

It is the following theorem, which presents the main result of the article.

\textbf{Theorem 3}
Let $\varphi_0\in D(A^{\beta+\frac{1}{2}})$ and $\varphi_1\in D(A^{\beta})$, where $\beta>\frac{s}{2}$, $M>0$ is a certain constant. Then, there are, depending on the value $M$, control $u(t,x)\in C([0,T]\times\Omega)$ and the time $T>0$, such that the solution of the problem (1)--(3) has the equalities
\begin{equation}
\label{12}
\theta(T,x)=\theta^{\prime}_t(T,x)=0,
\end{equation}
 and the restriction
$$
|u(t,x)|\leqslant M,
$$
for any $t\in (0,T]$, $x\in\Omega$ is done.

If we extend the constructed control function $u(t,x)$ by zero when $t>T$
then the system (1)--(3) will stays in the null state for $t>T$.

\begin{proof} Let $u(t,x)$ be the function, satisfying the theorem conditions, $T$ is some instant of time. According to the theorems 1 and 2, the solution of the task (\ref{1})--(\ref{3}) could be represented as (\ref{8}) and (\ref{9}). Hence we obtain
$$
\theta(t,x)=\frac{1}{\sqrt{2\pi}}\sum_{n=1}^{\infty}
\frac{(\varphi_{1n}+\lambda^+_n\varphi_{0n})e^{\lambda^+_nt}\psi_n(x)}{l^{(1)}_n(\lambda^+_n)}+
\frac{1}{\sqrt{2\pi}}\sum_{n=1}^{\infty}
\frac{(\varphi_{1n}+\lambda^-_n\varphi_{0n})e^{\lambda^-_nt}\psi_n(x)}{l^{(1)}_n(\lambda^-_n)}+
$$
$$
+\frac{1}{\sqrt{2\pi}}\sum_{n=1}^{\infty}\sum\limits_{k=0}^{N-1}
\left(\frac{(\varphi_{1n}-q_{k,n}\varphi_{0n})e^{-q_{k,n}t}}{l^{(1)}_n(-q_{k,n})}
\right)\psi_n(x)+
\frac{1}{\sqrt{2\pi}}\sum_{n=1}^{\infty}\frac{\int\limits_0^tu_n(s)e^{\lambda_n^+(t-s)}ds}
{l^{(1)}_n(\lambda_n^+)}\psi_n(x)+
$$
\begin{equation}
\label{13}
+\frac{1}{\sqrt{2\pi}}\sum_{n=1}^{\infty}\frac{\int\limits_0^tu_n(s)
e^{\lambda^-_n(t-s)}ds}{l^{(1)}_n(\lambda^-_n)}\psi_n(x)+
\frac{1}{\sqrt{2\pi}}\sum_{n=1}^{\infty}\sum\limits_{k=0}^{N-1}
\left(\frac{\int\limits_0^tu_n(s)e^{-q_{k,n}(t-s)}ds}{l^{(1)}_n(-q_{k,n})}
\right)\psi_n(x).
\end{equation}
Now we formally differentiate the last series by $t$ (uniform convergence of the series
for $\theta$ and $\theta_t$ will be proved in section 5):
$$
\frac{\partial\theta(t,x)}{\partial t}
=\frac{1}{\sqrt{2\pi}}\sum_{n=1}^{\infty}\frac{\lambda^+_n(\varphi_{1n}+\lambda^+_n\varphi_{0n})
e^{\lambda^+_nt}\psi_n(x)}{l^{(1)}_n(\lambda^+_n)}+
$$
$$
+\frac{1}{\sqrt{2\pi}}\sum_{n=1}^{\infty}
\frac{\lambda^-_n(\varphi_{1n}+\lambda^-_n\varphi_{0n})e^{\lambda^-_nt}\psi_n(x)}{l^{(1)}_n(\lambda^-_n)}+
$$
$$
+\frac{1}{\sqrt{2\pi}}\sum_{n=1}^{\infty}\sum\limits_{k=1}^{N-1}
\left(\frac{(-q_{k,n})(\varphi_{1n}-q_{k,n}\varphi_{0n})
e^{-q_{k,n}t}}{l^{(1)}_n(-q_{k,n})}
\right)\psi_n(x)+
$$
$$
+\frac{1}{\sqrt{2\pi}}\sum_{n=1}^{\infty}\left(\frac{1}{l_n^{(1)}(\lambda^+_n)}+
\frac{1}{l_n^{(1)}(\lambda^{-}_n)}+
\sum\limits_{k=0}^{N-1}\frac{1}{l_n^{(1)}(-q_{k,n})}\right)u_n(t)\psi_n(x)+
$$
$$
+\frac{1}{\sqrt{2\pi}}\sum_{n=1}^{\infty}\frac{\lambda^+_n\int\limits_0^tu_n(s)
e^{\lambda^+_n(t-s)}ds}{l^{(1)}_n(\lambda^+_n)}\psi_n(x)+
\frac{1}{\sqrt{2\pi}}\sum_{n=1}^{\infty}\frac{\lambda^-_n\int\limits_0^tu_n(s)
e^{\lambda^-_n(t-s)}ds}{l^{(1)}_n(\lambda^-_n)}\psi_n(x)+
$$
\begin{equation}
\label{14}
+\frac{1}{\sqrt{2\pi}}\sum_{n=1}^{\infty}\sum\limits_{k=1}^{N-1}
\left(\frac{(-q_{k,n})\int\limits_0^tu_n(s)e^{-q_{k,n}(t-s)}ds}{l^{(1)}_n(-q_{k,n})}
\right)\psi_n(x).
\end{equation}

Note that the fourth summand in (\ref{14}) is equal to zero, see lemma 1. Using conditions (\ref{12}) and formulae (\ref{13}), (\ref{14}) we obtain
$$
-\left(\frac{(\varphi_{1n}+\lambda^+_n\varphi_{0n})e^{\lambda^+_nT}}{l^{(1)}_n(\lambda^+_n)}+
\frac{(\varphi_{1n}+\lambda^-_n\varphi_{0n})e^{\lambda^-_nT}}{l^{(1)}_n(\lambda^-_n)}
+\sum\limits_{k=1}^{N-1}\frac{(\varphi_{1n}-q_{k,n}\varphi_{0n})e^{-q_{k,n}T}}{l^{(1)}_n(-q_{k,n})}\right)=
$$
$$
=\frac{\int\limits_0^Tu_n(s)e^{\lambda_n^+(T-s)}ds}{l^{(1)}_n(\lambda_n^+)}+
\frac{\int\limits_0^Tu_n(s)e^{\lambda^-_n(T-s)}ds}{l^{(1)}_n(\lambda^-_n)}+
$$
\begin{equation}
\label{15}
+\sum\limits_{k=0}^{N-1}\frac{\int\limits_0^Tu_n(s)e^{-q_{k,n}(T-s)}ds}{l^{(1)}_n(-q_{k,n})},\quad n=1,2,...,
\end{equation}
$$
-\frac{\lambda^+_n(\varphi_{1n}+\lambda^+_n\varphi_{0n})e^{\lambda^+_nT}}{l^{(1)}_n(\lambda^+_n)}-
\frac{\lambda^-_n(\varphi_{1n}+\lambda^-_n\varphi_{0n})e^{\lambda^-_nT}}{l^{(1)}_n(\lambda^-_n)}-
$$
$$
-\sum\limits_{k=1}^{N-1}\frac{(-q_{k,n})(\varphi_{1n}-q_{k,n}\varphi_{0n})e^{-q_{k,n}T}}{l^{(1)}_n(-q_{k,n})}=
$$
$$
=\frac{\lambda_n^+\int\limits_0^Tu_n(s)e^{\lambda_n^+(T-s)}ds}{l^{(1)}_n(\lambda_n^+)}+
\frac{\lambda^-_n\int\limits_0^Tu_n(s)e^{\lambda^-_n(T-s)}ds}{l^{(1)}_n(\lambda^-_n)}+
$$
\begin{equation}
\label{16}
+\sum\limits_{k=1}^{N-1}\frac{(-q_{k,n})\int\limits_0^Tu_n(s)e^{-q_{k,n}(T-s)}ds}{l^{(1)}_n(-q_{k,n})},\quad n=1,2,...\quad .
\end{equation}

We introduce
$$
a_n=-(\varphi_{1n}+\lambda_n^+\varphi_{0n}),\quad \bar{a}_n=-(\varphi_{1n}+\lambda_n^-\varphi_{0n}),\quad
$$
$$
b_{k,n}=-(\varphi_{1n}+(-q_{k,n})\varphi_{0n}),\quad k=0,1,2,...,N-1.
$$

Let equal coefficients preceding
$$
\frac{1}{l^{(1)}_n(\lambda_n^+)},\quad \frac{1}{l^{(1)}_n(\lambda_n^-)},\quad \frac{1}{l^{(1)}_n(-q_{k,n})},\quad k=0,1,2,...,N-1.
$$
in the right and left parts of equations (\ref{15}), (\ref{16}). Thus a new moments problem occurs:
$$
\int\limits_0^Tu_n(s)e^{\lambda_n^+(T-s)}ds=a_ne^{\lambda_n^+T},\quad \int\limits_0^Tu_n(s)e^{\lambda_n^-(T-s)}ds=\bar{a}_ne^{\lambda_n^-T},
\quad n=1,2,...,
$$
\begin{equation}
\label{17}
\int\limits_0^Tu_n(s)e^{-q_{k,n}(T-s)}ds=b_{k,n}e^{-q_{k,n}T},\quad k=0,1,2,...,N-1,\quad n=1,2,...\quad .
\end{equation}

Obviously, if moments problem (\ref{17}) is solvable then moments problem (\ref{15}), (\ref{16}) is solvable too. Elimination of
common factors in both parts (\ref{17})  allows us to represent this system as follows
$$
\int\limits_0^Tu_n(s)e^{-\lambda_n^+s}ds=a_n,\quad \int\limits_0^Tu_n(s)e^{-\lambda_n^-s}ds=\bar{a}_n,\quad n=1,2,...,
$$
\begin{equation}
\label{18}
\int\limits_0^Tu_n(s)e^{q_{k,n}s}ds=b_{k,n},\quad k=0,1,2,...,N-1,\quad n=1,2,...\quad .
\end{equation}

Let us replace $-\lambda_n^+=\lambda_n$ and $-\lambda_n^-=\overline{\lambda}_n$ in (\ref{18}). Notice that $Re\lambda_n>0$ and $q_{k,n}>0$, $k=1,2,...,N-1$ (see [6]). Finally, we obtain the system of $N+2$ moments for each natural number $n$:
$$
\int\limits_0^Tu_n(s)e^{\lambda_ns}ds=a_n,\quad \int\limits_0^Tu_n(s)e^{\overline{\lambda}_ns}ds=\bar{a}_n,\quad n=1,2,...,
$$
\begin{equation}
\label{20}
\int\limits_0^Tu_n(s)e^{q_{k,n}s}ds=b_{k,n},\quad k=0,1,2,...,N-1,\quad n=1,2,...\quad .
\end{equation}

The solution of (\ref{20}) is sought in the following form:
\begin{equation}
\label{21}
u_n(s)=C_{-2,n}e^{\lambda_ns}+C_{-1,n}e^{\overline{\lambda}_ns}+\sum\limits_{j=0}^{N-1}C_{j,n}e^{q_{j,n}s},\quad n=1,2,...\quad .
\end{equation}
Set $C_{-2,n}$ $C_{-1,n}$ and $C_{k,n}$ as some unknown constants. Substituting (\ref{21}) in (\ref{20}),
we get the system of $N+2$ algebraic equations for each natural number $n$:
$$
C_{-2,n}\int\limits_0^Te^{2\lambda_ns}ds+C_{-1,n}\int\limits_0^Te^{(\lambda_n+\overline{\lambda}_n)s}ds+
\sum\limits_{k=0}^{N-1}C_{k,n}\int\limits_0^Te^{(\lambda_n+q_{k,n})s}ds=a_n,
$$
$$
C_{-2,n}\int\limits_0^Te^{(\overline{\lambda}_n+\lambda_n)s}ds+C_{-1,n}\int\limits_0^Te^{2\overline{\lambda}_ns}ds+
\sum\limits_{k=0}^{N-1}C_{k,n}\int\limits_0^Te^{(\overline{\lambda}_n+q_{k,n})s}ds=a_n,
$$
\begin{equation}
\label{22}
C_{-2,n}\int\limits_0^Te^{(\lambda_n+q_{k,n})s}ds+C_{-1,n}\int\limits_0^Te^{(\overline{\lambda}_n+q_{k,n})s}ds
+\sum\limits_{j=0}^{N-1}C_{j,n}\int\limits_0^Te^{(q_{j,n}+q_{k,n})s}ds=b_{k,n},\quad k=0,1,2,...,N-1.
\end{equation}

Let us find the determinant $\Delta_n$ of the problem (\ref{22}).
$$
\small
\Delta_n=
\left|\begin{array}{cccccc}
\int\limits_0^Te^{2\lambda_ns}ds& \int\limits_0^Te^{(\lambda_n+\overline{\lambda}_n)s}ds&
\int\limits_0^Te^{\lambda_n s}ds&
\int\limits_0^Te^{(\lambda_n+q_{1,n})s}ds&\ldots & \int\limits_0^Te^{(\lambda_n+q_{N-1,n})s}ds\\
\int\limits_0^Te^{(\overline{\lambda}_n+\lambda_n)s}ds& \int\limits_0^Te^{2\overline{\lambda}_ns}ds&
\int\limits_0^Te^{\overline{\lambda}_n s}ds&
\int\limits_0^Te^{(\overline{\lambda}_n+q_{1,n})s}ds&\ldots &
\int\limits_0^Te^{(\overline{\lambda}_n+q_{N-1,n})s}ds\\
\int\limits_0^Te^{\lambda_n s}ds &
\int\limits_0^Te^{\overline{\lambda}_n s}ds &
T & \int\limits_0^Te^{q_{1,n}s}ds &
\ldots & \int\limits_0^Te^{q_{N-1,n} s}ds
\\
\int\limits_0^Te^{(q_{1,n}+\lambda_n)s}ds& \int\limits_0^Te^{(q_{1,n}+\overline{\lambda}_n)s}ds&
\int\limits_0^Te^{q_{1,n} s}ds&
\int\limits_0^Te^{2q_{1,n}s}ds&\ldots & \int\limits_0^Te^{(q_{1,n}+q_{N-1,n})s}ds
\\
\vdots& \vdots& \vdots& \vdots& \ddots& \vdots\\
\int\limits_0^Te^{(q_{N-1,n}+\lambda_n)s}ds& \int\limits_0^Te^{(q_{N-1,n}+\overline{\lambda}_n)s}ds&
\int\limits_0^Te^{q_{N-1,n} s}ds&
\int\limits_0^Te^{(q_{N-1,n}+q_{1,n})s}ds&\ldots & \int\limits_0^Te^{2q_{N-1,n}s}ds
\end{array}\right|
$$

We notice that all determinants $\Delta_n$ are nonzero for any natural index $n$ because
$\Delta_n$ is the Gram determinant.

As far as
\begin{equation}
\label{23}
\int\limits_0^Te^{(q_{i,n}+q_{j,n})s}ds=\frac{1}{q_{i,n}+q_{j,n}}e^{(q_{i,n}+q_{j,n})T}-\frac{1}{q_{i,n}+q_{j,n}},
\end{equation}
then using the equality (\ref{23}) and a well known property of determinants
\begin{gather*}
\small
\left|\begin{array}{cccc}
a_{11}& a_{12}&\ldots & a_{1n} \\
a_{21}& a_{22}&\ldots & a_{2n} \\
\vdots& \vdots&  &\vdots\\
b_{i1}+c_{i1}& b_{i2}+c_{i2}&\ldots & b_{in}+c_{in} \\
\vdots& \vdots& \ddots& \vdots\\
a_{n1}& a_{n2}&\ldots & a_{nn}
\end{array}\right|=
\small
\left|\begin{array}{cccc}
a_{11}& a_{12}&\ldots & a_{1n} \\
a_{21}& a_{22}&\ldots & a_{2n} \\
\vdots& \vdots&  &\vdots\\
b_{i1}& b_{i2}&\ldots & b_{in} \\
\vdots& \vdots& \ddots& \vdots\\
a_{n1}& a_{n2}&\ldots & a_{nn}
\end{array}\right|+
\end{gather*}
\begin{gather}
\small
\label{24}
+\left|\begin{array}{cccc}
a_{11}& a_{12}&\ldots & a_{1n} \\
a_{21}& a_{22}&\ldots & a_{2n} \\
\vdots& \vdots&  &\vdots\\
c_{i1}& c_{i2}&\ldots & c_{in} \\
\vdots& \vdots& \ddots& \vdots\\
a_{n1}& a_{n2}&\ldots & a_{nn}
\end{array}\right|,
\end{gather}
we obtain that
\begin{gather}
\small
\label{25}
\Delta_n=\left|\begin{array}{cccccc}
\frac{e^{2\lambda_nT}}{2\lambda_n}& \frac{e^{(\lambda_n+\overline{\lambda}_n)T}}{\lambda_n+\overline{\lambda}_n}&
\frac{e^{\lambda_n T}}{\lambda_n} &
\frac{e^{(\lambda_n+q_{1,n})T}}{\lambda_n+q_{1,n}}&\ldots &
\frac{e^{(\lambda_n+q_{N-1,n})T}}{\lambda_n+q_{N-1,n}}\\
\\
\frac{e^{(\overline{\lambda}_n+\lambda_n)T}}{\overline{\lambda}_n+\lambda_n}& \frac{e^{2\overline{\lambda}_n T}}{2\overline{\lambda}_n}&
\frac{e^{\overline{\lambda}_n T}}{\overline{\lambda}_n}&
\frac{e^{(\overline{\lambda}_n+q_{1,n})T}}{\overline{\lambda}_n+q_{1,n}}&\ldots &
\frac{e^{(\overline{\lambda}_n+q_{N-1,n})T}}{\overline{\lambda}_n+q_{N-1,n}}\\
\\
\frac{e^{\lambda_n T}}{\lambda_n}& \frac{e^{\overline{\lambda}_n T}}{\overline{\lambda}_n}& T & \frac{e^{q_{1,n}T}}{q_{1,n}}&\ldots &
\frac{e^{q_{N-1,n}T}}{q_{N-1,n}}\\
\\
\frac{e^{(q_{1,n}+\lambda_n)T}}{q_{1,n}+\lambda_n}& \frac{e^{(q_{1,n}+\overline{\lambda}_n)T}}{q_{1,n}+\overline{\lambda}_n}& \frac{e^{q_{1,n}T}}{q_{1,n}}&
\frac{e^{2q_{1,n}T}}{2q_{1,n}}&
\ldots &
\frac{e^{(q_{1,n}+q_{N-1,n})T}}{q_{1,n}+q_{N-1,n}}\\
\vdots& \vdots& \vdots& \vdots& \ddots& \vdots\\
\frac{e^{(q_{N-1,n}+\lambda_n)T}}{q_{N-1,n}+\lambda_n}& \frac{e^{(q_{N-1,n}+\overline{\lambda}_n)T}}{q_{N-1,n}+\overline{\lambda}_n}&
\frac{e^{q_{N-1,n} T}}{q_{N-1,n}}&
\frac{e^{(q_{N-1,n}+q_{1,n})T}}{q_{N-1,n}+q_{1,n}}& \ldots& \frac{e^{2q_{N-1,n}T}}{2q_{N-1,n}}
\end{array}\right|+\beta_n(T),
\end{gather}
 where $\beta_n(T)$ is a sum of all other determinants, which are the result of $N+2$-fold application of the property (\ref{24}) to each row of the determinant $\Delta_n$.

Let us factor out $e^{\lambda_nT}$ from the first row of the determinant in the right part of (\ref{25}), and then take out the same factor from the first column, now the similar action can be made for the second row and column with the factor $e^{\overline{\lambda}_nT}$, and so on.

Thus we get that

\begin{gather}
\small
\label{26}
\Delta_n=e^{2\lambda_nT}e^{2\overline{\lambda}_n T}\prod\limits_{j=1}^{N-1}e^{2q_{j,n}T}
\left|\begin{array}{cccccc}
\frac{1}{2\lambda_n}& \frac{1}{\lambda_n+\overline{\lambda}_n}&
\frac{1}{\lambda_n} &
\frac{1}{\lambda_n+q_{1,n}}&\ldots &
\frac{1}{\lambda_n+q_{N-1,n}}\\
\\
\frac{1}{\overline{\lambda}_n+\lambda_n}& \frac{1}{2\overline{\lambda}_n}&
\frac{1}{\overline{\lambda}_n}&
\frac{1}{\overline{\lambda}_n+q_{1,n}}&\ldots &
\frac{1}{\overline{\lambda}_n+q_{N-1,n}}\\
\\
\frac{1}{\lambda_n}& \frac{1}{\overline{\lambda}_n}& T & \frac{1}{q_{1,n}}&\ldots &
\frac{1}{q_{N-1,n}}\\
\\
\frac{1}{q_{1,n}+\lambda_n}& \frac{1}{q_{1,n}+\overline{\lambda}_n}& \frac{1}{q_{1,n}}&
\frac{1}{2q_{1,n}}&
\ldots &
\frac{1}{q_{1,n}+q_{N-1,n}}\\
\vdots& \vdots& \vdots& \vdots& \ddots& \vdots\\
\frac{1}{q_{N-1,n}+\lambda_n}& \frac{1}{q_{N-1,n}+\overline{\lambda}_n}&
\frac{1}{q_{N-1,n}}&
\frac{1}{q_{N-1,n}+q_{1,n}}& \ldots& \frac{1}{2q_{N-1,n}}
\end{array}\right|
+\beta_n(T).
\end{gather}

Denote
$$
\bar{\Delta}_n=
\small
\left|\begin{array}{cccccc}
\frac{1}{2\lambda_n}& \frac{1}{\lambda_n+\overline{\lambda}_n}&
\frac{1}{\lambda_n} &
\frac{1}{\lambda_n+q_{1,n}}&\ldots &
\frac{1}{\lambda_n+q_{N-1,n}}\\
\\
\frac{1}{\overline{\lambda}_n+\lambda_n}& \frac{1}{2\overline{\lambda}_n}&
\frac{1}{\overline{\lambda}_n}&
\frac{1}{\overline{\lambda}_n+q_{1,n}}&\ldots &
\frac{1}{\overline{\lambda}_n+q_{N-1,n}}\\
\\
\frac{1}{\lambda_n}& \frac{1}{\overline{\lambda}_n}& T & \frac{1}{q_{1,n}}&\ldots &
\frac{1}{q_{N-1,n}}\\
\\
\frac{1}{q_{1,n}+\lambda_n}& \frac{1}{q_{1,n}+\overline{\lambda}_n}& \frac{1}{q_{1,n}}&
\frac{1}{2q_{1,n}}&
\ldots &
\frac{1}{q_{1,n}+q_{N-1,n}}\\
\vdots& \vdots& \vdots& \vdots& \ddots& \vdots\\
\frac{1}{q_{N-1,n}+\lambda_n}& \frac{1}{q_{N-1,n}+\overline{\lambda}_n}&
\frac{1}{q_{N-1,n}}&
\frac{1}{q_{N-1,n}+q_{1,n}}& \ldots& \frac{1}{2q_{N-1,n}}
\end{array}\right|.
$$

Then
$$
\Delta_n=e^{2\lambda_nT}e^{2\overline{\lambda}_nT}
\prod\limits_{j=1}^{N-1}e^{2q_{j,n}T}\left(\bar{\Delta}_n+
e^{-2\lambda_nT}e^{-2\overline{\lambda}_nT}\prod\limits_{j=1}^{N-1}e^{-2q_{j,n}T}
\beta_n(T)\right).
$$

We notice that the sequence of the modules of complex roots $\{|\lambda_n|\}$ tends to $+\infty$ if $n\rightarrow+\infty$ but $Re\lambda_n=\mu+O(n^{-2})$ ($\mu>0$) and
the sequence of real numbers $\{q_{k,n}\}_{n=1}^{\infty}$ converges to some positive number $q_k$, actually $q_{k,n}=q_k+O(n^{-2})$ (see \cite{Ivanov_Sheronova}).
Then in virtue of the definition of $\beta_n(T)$, the following fact takes place:
$$
\left|e^{-2\lambda_nT}e^{-2\overline{\lambda}_nT}\prod\limits_{j=1}^{N-1}e^{-2q_{j,n}T}\beta_n(T)\right|\rightarrow 0,\:\: T\rightarrow +\infty.
$$

Let us represent the determinant $\bar{\Delta}_n$ in the following form

\begin{gather}
\small
\label{27}
\bar{\Delta}_n=
\small
\frac{1}{(2Re\lambda_n)^2}\left|\begin{array}{cccc}
T & \frac{1}{q_{1,n}}&\ldots &
\frac{1}{q_{N-1,n}}\\
\frac{1}{q_{1,n}}& \frac{1}{2q_{1,n}}&
\ldots &
\frac{1}{q_{1,n}+q_{N-1,n}}\\
\vdots& \vdots& \ddots& \vdots\\
\frac{1}{q_{N-1,n}}&
\frac{1}{q_{N-1,n}+q_{1,n}}& \ldots& \frac{1}{2q_{N-1,n}}
\end{array}\right|+\Lambda_n(T),
\end{gather}
where $\Lambda_n(T)\rightarrow 0$, $n\rightarrow+\infty$ when $T$ is fixed.

Let us make the following notation
$$
P_n=
\small
\left|\begin{array}{cccc}
\frac{1}{2q_{1,n}}& \frac{1}{q_{1,n}+q_{2,n}}&\ldots &
\frac{1}{q_{1,n}+q_{N-1,n}}\\
\frac{1}{q_{2,n}+q_{1,n}}& \frac{1}{2q_{2,n}}&\ldots &
\frac{1}{q_{2,n}+q_{N-1,n}}\\
\vdots& \vdots& \ddots& \vdots\\
\frac{1}{q_{N-1,n}+q_{1,n}}& \frac{1}{q_{N-1,n}+q_{2,n}}&\ldots& \frac{1}{2q_{N-1,n}}
\end{array}\right|.
$$

$P_n$ is the Cauchy determinant. It is a well known fact that
$$
P_n=\frac{\prod\limits_{N-1\geqslant i>j\geqslant 1}(q_{i,n}-q_{j,n})^2}{\prod\limits_{i,j=1}^{N-1}(q_{i,n}+q_{j,n})}.
$$

As far as  $q_{i,n}$, $i=1,2,...,N-1$ are pairwise different for any $n$ (see [6]), then  $P_n$ is nonzero.
It is obvious that
$$
\left|\begin{array}{cccc}
T & \frac{1}{q_{1,n}}&\ldots &
\frac{1}{q_{N-1,n}}\\
\frac{1}{q_{1,n}}& \frac{1}{2q_{1,n}}&
\ldots &
\frac{1}{q_{1,n}+q_{N-1,n}}\\
\vdots& \vdots& \ddots& \vdots\\
\frac{1}{q_{N-1,n}}&
\frac{1}{q_{N-1,n}+q_{1,n}}& \ldots& \frac{1}{2q_{N-1,n}}
\end{array}\right|=
T\left(P_n+\xi_n(T)\right),
$$
where $\xi_n(T)\rightarrow 0$, $T\rightarrow+\infty$.

Hence,
$$
\Delta_n=e^{2\lambda_nT}e^{2\overline{\lambda}_nT}
\prod\limits_{j=1}^{N-1}e^{2q_{j,n}T}
\left(\frac{1}{(2Re\lambda_n)^2}T\left(P_n+\xi_n(T)\right)
+\Lambda_n(T)+e^{-2\lambda_nT}e^{-2\overline{\lambda}_nT}
\prod\limits_{j=1}^{N-1}e^{-2q_{j,n}T}
\beta_n(T)\right)=
$$
$$
\frac{TP_n}{(2Re\lambda_n)^2}e^{2(\lambda_n+\overline{\lambda}_n)T}
\prod\limits_{j=1}^{N-1}e^{2q_{j,n}T}\left(1+\frac{\xi_n(T)}{P_n}+
\frac{(2Re\lambda_n)^2}{TP_n}\Lambda_n(T)
+\frac{(2Re\lambda_n)^2}{TP_n}e^{-2(\lambda_n+\overline{\lambda}_n)T}\prod\limits_{j=1}^{N-1}e^{-2q_{j,n}T}
\beta_n(T)\right).
$$

Let us denote
$$
\bar{\xi}_n(T)=\frac{\xi_n(T)}{P_n},\quad
\bar{\Lambda}_n(T)=\frac{(2Re\lambda_n)^2}{TP_n}\Lambda_n(T),
$$
$$
\bar{\beta}_n(T)=
\frac{(2Re\lambda_n)^2}{TP_n}e^{-2(\lambda_n+\overline{\lambda}_n)T}
\prod\limits_{j=1}^{N-1}e^{-2q_{j,n}T}\beta_n(T).
$$

It leads to the following:
\begin{equation}
\label{28}
\Delta_n=\frac{TP_n}{(2Re\lambda_n)^2}e^{2\lambda_nT}e^{2\overline{\lambda}_nT}
\prod\limits_{j=1}^{N-1}e^{2q_{j,n}T}\left(1+\bar{\xi}_n(T)+\bar{\Lambda}_n(T)
+\bar{\beta}_n(T)\right).
\end{equation}

Notice that $\bar{\Lambda}_n(T)\rightarrow 0$ when $n\rightarrow+\infty$ (uniformly by $T\in[T_{\ast},+\infty)$
for any $T_{\ast}>0$) and $\bar{\xi}_n(T)$,
$\bar{\beta}_n(T)\rightarrow 0$ if $T\rightarrow+\infty$, more exactly for any $\varepsilon>0$ there exist $T>0$ and $n_{\ast}$ such that
$$
|\bar{\xi}_n(T)|<\varepsilon,\:\:|\bar{\beta}_n(T)|<\varepsilon
$$
for any $n>n_{\ast}$.

Let us determine $\Delta_{-2,n}$ by the formula:
$$
\small
\Delta_{-2,n}=
\left|\begin{array}{cccccc}
a_n& \int\limits_0^Te^{(\lambda_n+\overline{\lambda}_n)s}ds&
\int\limits_0^Te^{\lambda_n s}ds&
\int\limits_0^Te^{(\lambda_n+q_{1,n})s}ds&\ldots & \int\limits_0^Te^{(\lambda_n+q_{N-1,n})s}ds\\
\overline{a}_n& \int\limits_0^Te^{2\overline{\lambda}_ns}ds&
\int\limits_0^Te^{\overline{\lambda}_n s}ds&
\int\limits_0^Te^{(\overline{\lambda}_n+q_{1,n})s}ds&\ldots &
\int\limits_0^Te^{(\overline{\lambda}_n+q_{N-1,n})s}ds\\
b_{0,n} &
\int\limits_0^Te^{\overline{\lambda}_n s}ds &
T & \int\limits_0^Te^{q_{1,n}s}ds &
\ldots & \int\limits_0^Te^{q_{N-1,n} s}ds
\\
b_{1,n}& \int\limits_0^Te^{(q_{1,n}+\overline{\lambda}_n)s}ds&
\int\limits_0^Te^{q_{1,n} s}ds&
\int\limits_0^Te^{2q_{1,n}s}ds&\ldots & \int\limits_0^Te^{(q_{1,n}+q_{N-1,n})s}ds
\\
\vdots& \vdots& \vdots& \vdots& \ddots& \vdots\\
b_{N-1,n}& \int\limits_0^Te^{(q_{N-1,n}+\overline{\lambda}_n)s}ds&
\int\limits_0^Te^{q_{N-1,n} s}ds&
\int\limits_0^Te^{(q_{N-1,n}+q_{1,n})s}ds&\ldots & \int\limits_0^Te^{2q_{N-1,n}s}ds
\end{array}\right|.
$$

Set likewise $\Delta_{k,n}$, where $k=-1,0,1,2,...,N-1$:
$$
\small
\Delta_{k,n}=
\left|\begin{array}{ccccc}
\int\limits_0^Te^{2\lambda_ns}ds& \ldots & a_n& \ldots & \int\limits_0^Te^{(\lambda_n+q_{N-1,n})s}ds\\
\int\limits_0^Te^{(\overline{\lambda}_n+\lambda_n)s}ds& \ldots & \overline{a}_n & \ldots &
\int\limits_0^Te^{(\overline{\lambda}_n+q_{N-1,n})s}ds\\
\int\limits_0^Te^{(q_{0,n}+\lambda_n)s}ds& \ldots & b_{0,n} & \ldots &
\int\limits_0^Te^{(q_{0,n}+q_{N-1,n})s}ds\\
\vdots& \vdots& \vdots&  \ddots& \vdots\\
\int\limits_0^Te^{(q_{N-1,n}+\lambda_n)s}ds& \ldots & b_{N-1,n} & \ldots &
\int\limits_0^Te^{2q_{N-1,n}s}ds
\end{array}\right|,
$$
where the column $\{a_n,\overline{a}_n,b_{0,n},b_{1,n},...,b_{N-1,n}\}$ takes the $k$-th place.

Applying Cramer's rule, we obtain:
$$
C_{-2,n}=\frac{\Delta_{-2,n}}{\Delta_n},\quad C_{-1,n}=\frac{\Delta_{-1,n}}{\Delta_n}, \quad C_{k,n}=\frac{\Delta_{k,n}}{\Delta_n},\quad k=0,1,2,...,N-1.
$$

Thus the solution of (\ref{20}) in the instant of time $t$ has the following form:
$$
u_n(t)=\frac{\Delta_{-2,n}}{\Delta_n}e^{\lambda_nt}+\frac{\Delta_{-1,n}}{\Delta_n}e^{\overline{\lambda}_nt}+
\sum\limits_{k=0}^{N-1}\frac{\Delta_{k,n}}{\Delta_n}e^{q_{k,n}t},
$$

Let $\lambda_n=\mu_n-i\nu_n$. The article [6] proves, that $\mu_n,\nu_n>0$ for any natural index $n$.
The estimation of the modulus of the function $u_n(t)$ for any natural $n$ should be provided. So we have
\begin{equation}
\label{29}
|u_n(t)|\leqslant
\frac{|\Delta_{-2,n}|}{|\Delta_n|}e^{\mu_nT}+\frac{|\Delta_{-1,n}|}{|\Delta_n|}e^{\mu_nT}+
\sum\limits_{k=0}^{N-1}\frac{|\Delta_{k,n}|}{|\Delta_n|}e^{q_{k,n}T}.
\end{equation}

Calculating determinants $\Delta_{-2,n}$, $\Delta_{-1,n}$, $\Delta_{k,n}$, $k=0,1,2,...,N-1$, it is clear that the part of summands consists of the different exponential product. Notice that the exponential product with the largest number of exponential factors in the determinant $\Delta_{-2,n}$ has the form
$$
e^{\lambda_nT}e^{2\overline{\lambda}_nT}e^{2q_{1,n}T}e^{2q_{2,n}T}
\cdot\cdot\cdot e^{2q_{N-1,n}T}.
$$
Similarly for $\Delta_{-1,n}$. For $\Delta_{k,n}$ ($k\neq 0$) we have:
$$
e^{2\lambda_nT}e^{2\overline{\lambda}_nT}e^{2q_{1,n}T}e^{2q_{2,n}T}
\cdot\cdot\cdot e^{q_{k,n}T}\cdot\cdot\cdot e^{2q_{N-1,n}T}.
$$
And in $\Delta_{0,n}$:
$$
e^{2\lambda_nT}e^{2\overline{\lambda}_nT}e^{2q_{1,n}T}e^{2q_{2,n}T}
\cdot\cdot\cdot e^{2q_{N-1,n}T}.
$$
It means that $u_n$ decreases as $T^{-1}$ when $T\rightarrow+\infty$.
Thus it is possible to make the modulus of the function $u_n(t)$, and hence of the control $u(t)$, be indefinitely small by means of increasing time control. Using (\ref{28}), (\ref{29}) we obtain

$$
|u_n(t)|\leqslant
\frac{4\mu^2_n}{T|P_n|e^{4\mu_nT}\prod\limits_{j=1}^{N-1}e^{2q_{j,n}T}\left(1-|\bar{\xi}_n(T)|
-|\bar{\Lambda}_n(T)|
-|\bar{\beta}_n(T)|\right)}e^{\mu_nT}(|\Delta_{-2,n}|+|\Delta_{-1,n}|)+
$$
\begin{equation}
\label{30}
+\sum\limits_{k=0}^{N-1}\frac{4\mu^2_n|\Delta_{k,n}|}{T|P_n|e^{4\mu_nT}
\prod\limits_{j=1}^{N-1}e^{2q_{j,n}T}\left(1-|\bar{\xi}_n(T)|-|\bar{\Lambda}_n(T)|
-|\bar{\beta}_n(T)|\right)}e^{q_{k,n}T},\quad t\in[0,T].
\end{equation}

Strictly speaking, evaluation (\ref{30}) is obtained for numbers $n$ greater than some number $n_{\ast}$.
Obviously first members of the series are estimated from above on the module by constant $c_{\ast}/T$.

Using (\ref{29}) and (\ref{30}) let us prove that there is the time required to stabilise the system, providing that the function $u(t,x)$ satisfies the condition
\begin{equation}
\label{32}
|u(t,x)|\leqslant M,
\end{equation}
where $M$ is an arbitrary constant.

As far as the sequences of real numbers  $\{\mu_n\}$, $\{\nu_n\}$, $\{q_{k,n}\}$ are such that $\mu_n=\mu+O(n^{-2})$, $\nu_n=D\alpha_n$ and $q_{k,n}=q_k+O(n^{-2})$,
where $\mu$, $D$, $q_k$ are some positive numbers (see [6]),
and moreover, the sequences $\{|a_n|\}$, $\{|b_{k,n}|\}$, $\{|\Lambda_n|\}$
tend to zero, then the following estimate takes place:
\begin{equation}
\label{33}
|u(t,x)|\leqslant\frac{c}{T}\sqrt{\sum_{n=1}^{\infty}\alpha_n^{2\beta}\left(|a_n|^2+|\overline{a}_n|^2+
\sum_{k=0}^{N-1}|b_{k,n}|^2\right)}\sqrt{\sum_{n=1}^{\infty}\alpha_n^{-2\beta}\psi_n^2(x)},
\end{equation}
 where $c$ is some constant and $T$ is large enough.
It is a well known fact that (see \cite{Agmon})
$$
\sum_{n=1}^{\infty}\alpha_n^{-2\beta}\psi_n^2(x)\leqslant\mathrm{const}\quad\mbox{if}\:\: 2\beta>s.
$$
Moreover, the series $\sum\limits_{n=1}^{\infty}\alpha_n^{-2\beta}\psi_n^2(x)$ is a continuous
function (if $2\beta>s$) and
$$
\sum_{n=1}^{\infty}\alpha_n^{2\beta+2}|\varphi_{0n}|^2=
\int\limits_{\Omega}\left|A^{\frac{\beta+1}{2}}\varphi_{0}(x)\right|^2dx,\quad
\sum_{n=1}^{\infty}\alpha_n^{2\beta}|\varphi_{1n}|^2=
\int\limits_{\Omega}\left|A^{\frac{\beta}{2}}\varphi_{1}(x)\right|^2dx.
$$
The latter series converges if $A^{\frac{\beta+1}{2}}\varphi_{0}(x)\in L_2(\Omega)$ and $A^{\frac{\beta}{2}}\varphi_{1}(x)\in L_2(\Omega)$ then it must be $\varphi_{0}(x)\in D(A^{\frac{\beta+1}{2}})$ and
$\varphi_{1}(x)\in D(A^{\frac{\beta}{2}})$.
But in the theorem statement these conditions were imposed on the initial data
as $D(A^{\beta+\frac{1}{2}})\subset D(A^{\frac{\beta+1}{2}})$ and $D(A^{\beta})\subset D(A^{\frac{\beta}{2}})$.

Thus we obtain
\begin{equation}
\label{34}
|u(t,x)|\leqslant\frac{C_1}{T}\leqslant M,
\end{equation}
where $c_1$ is some constant and $T$ is great enough.

Let us show that the following statement takes place: $u(t,x)\in C([0,T]\times\Omega)$. We obtain
\begin{equation}
\label{31}
|u(t,x)|\leqslant \sum_{n=1}^{\infty}|u_n(t)||\psi_n(x)|
\leqslant\frac{1}{2}\left(\sum_{n=1}^{\infty}\alpha_n^{2\beta}|u_n(t)|^2+
\sum_{n=1}^{\infty}\alpha_n^{-2\beta}\psi_n^2(x)
\right)
\end{equation}
Note that series $\sum\limits_{n=1}^{\infty}\alpha_n^{-2\beta}\psi_n^2(x)$
converges uniformly, it follows from the Dini's theorem.
Then, using the sign of Weierstrass, we obtain the required statement of the theorem.

\end{proof}

\section{Controllability to rest}

Let us show that the control function, constructed in the proof of the theorem, in fact, drives the system to rest.
For this purpose we use the formula (\ref{13}) and integral equations (\ref{18}).
The function $u(t,x)$ can be continued by zero at $t>T$, then the system (\ref{18})
can be written (when $t>T$) as:
$$
\int\limits_0^{t}u_n(s)e^{-\lambda_n^+s}ds=a_n,\quad
\int\limits_0^{t}u_n(s)e^{-\lambda_n^-s}ds=\bar{a}_n,\quad n=1,2,...,
$$
\begin{equation}
\label{35}
\int\limits_0^{t}u_n(s)e^{q_{k,n}s}ds=b_{k,n},\quad k=0,1,2,...,N-1,\quad n=1,2,...\quad .
\end{equation}
It means that values of the Laplace transform of function $u_n(t)$ at nulls of function $l_n(\lambda)$ are equal to the given numbers. Further, taking $t>T$ in formula (\ref{13}) and using (\ref{35}),we obtain that $\theta(t,x)\equiv 0$, at any $t>T$.

\textbf{Remark.}
The article \cite{Romanov}, devoted to the analogical problem (\ref{1})--(\ref{3})
(for one-dimensional system), did not take into account a zero root of $l_n(\lambda)$.
That is why in \cite{Romanov} the right part of the estimate analogous to (\ref{34}) decreases exponentially and not as $T^{-1}$. In general, it is not true. Exponential decreasing, in the given problem, takes place when the second initial condition is equal to zero.

\section{Existence of the solution and convergence of series}

Now we show that $\theta(t,x)$ can be understood as a solution of problem the
(\ref{1})--(\ref{3}) at $t>T$. For this purpose let us use the results in
\cite{Vlasov_1}. According to this article, for the existence and uniqueness of a solution of the problem (\ref{1})--(\ref{3}) at $t\in[0,+\infty)$ the following conditions of smoothness should take place: $\varphi_0\in D(A)$, $\varphi_1\in D(A^{\frac{1}{2}})$
and $A^{\frac{1}{2}}u(t,x)\in L_2(R_+,H)$. In this case, solution $\theta(t,x)$
is in space $W^2_{2,\gamma}(R_+,A)$ for any $\gamma>0$.

The conditions of smoothness are satisfied in virtue of the statement of theorem 3.
Let us show now that the right part is smooth enough. We have:
$$
\int\limits_{\Omega}\left|A^{\frac{1}{2}}u(t,x)\right|^2dx=
\sum_{n=1}^{\infty}\alpha^2_n |u_{n}(t)|^2.
$$
Using evaluations for $|u_n(t)|^2$ which are analogous to (\ref{33}), we get that
$$
\int\limits_{\Omega}\left|A^{\frac{1}{2}}u(t,x)\right|^2dx\leqslant C_2,\quad t\in[0,T]
$$
because
$$
\sum_{n=1}^{\infty}\alpha_n^{4\beta+2}|\varphi_{0n}|^2=
\int\limits_{\Omega}\left|A^{\beta+\frac{1}{2}}\varphi_{0}(x)\right|^2dx,\quad
\sum_{n=1}^{\infty}\alpha_n^{4\beta}|\varphi_{1n}|^2=
\int\limits_{\Omega}\left|A^{\beta}\varphi_{1}(x)\right|^2dx
$$
and $4\beta>2s\geqslant 2$.
As $u(t,x)$ is equal to zero at $t>T$, then $A^{\frac{1}{2}}u(t,x)\in L_2(R_+,H)$.

Now we show the uniform convergence with respect to $t\in[0,T]$ in the norm of $L_2(\Omega)$
of the series for $\theta_t$. For solution $\theta$ the proof
is analogous and easier. We get:
$$
l_n(\lambda)=\lambda^2+\alpha_n^2\sum\limits_{k=1}^{N}\frac{c_k}{\gamma_k}
-\alpha_n^2\sum\limits_{k=1}^{N}\frac{c_k}{\gamma_k(\lambda+\gamma_k)}.
$$
Then the first derivative has the form:
$$
l^{(1)}_n(\lambda)=2\lambda
+\alpha_n^2\sum\limits_{k=1}^{N}\frac{c_k}{\gamma_k(\lambda+\gamma_k)^2}.
$$
Consider the series (\ref{14}), $t\in[0,T]$. Using Parseval's identity, we get:
$$
\sum_{n=1}^{\infty}\left|\frac{\lambda^+_n(\varphi_{1n}+\lambda^+_n\varphi_{0n})
e^{\lambda^+_nt}}{l^{(1)}_n(\lambda^+_n)}\right|^2\leqslant
C_3\sum_{n=1}^{\infty}\frac{\alpha^2_n(|\varphi_{1n}|^2+\alpha^2_n|\varphi_{0n}|^2)}
{\alpha_n^2}=
C_3\sum_{n=1}^{\infty}(|\varphi_{1n}|^2+\alpha^2_n|\varphi_{0n}|^2).
$$
The last numerical series converges due to the choice of initial data spaces. Likewise,
the proof is carried out for $\lambda=\lambda^-_n$.
Further:
$$
\sum_{n=1}^{\infty}\left|\sum\limits_{k=1}^{N-1}
\frac{(-q_{k,n})(\varphi_{1n}-q_{k,n}\varphi_{0n})
e^{-q_{k,n}t}}{l^{(1)}_n(-q_{k,n})}\right|^2\leqslant C_4
\sum_{n=1}^{\infty}\left(\sum\limits_{k=1}^{N-1}\left|
\frac{(-q_{k,n})(\varphi_{1n}-q_{k,n}\varphi_{0n})
e^{-q_{k,n}t}}{l^{(1)}_n(-q_{k,n})}\right|\right)^2\leqslant
$$
$$
C_5\sum_{n=1}^{\infty}(|\varphi_{1n}|^2+|\varphi_{0n}|^2).
$$
The last series converges obviously.

$$
\sum_{n=1}^{\infty}\left|\frac{\lambda^+_n\int\limits_0^tu_n(s)
e^{\lambda^+_n(t-s)}ds}{l^{(1)}_n(\lambda^+_n)}\right|^2\leqslant
C_6\sum_{n=1}^{\infty}\frac{\alpha^2_n\int\limits_0^t|u_n(s)|^2ds
\int\limits_0^t|e^{\lambda^+_n(t-s)}|^2ds}{\alpha^2_n}\leqslant
C_7\sum_{n=1}^{\infty}\int\limits_0^T|u_n(s)|^2ds.
$$
The last numerical series converges, because it was proved earlier that
$u(t,x)\in C([0,T]\times\Omega)$.
Analogously for $\lambda=\lambda^-_n$.

$$
\sum_{n=1}^{\infty}\left|\sum\limits_{k=1}^{N-1}
\frac{(-q_{k,n})\int\limits_0^tu_n(s)e^{-q_{k,n}(t-s)}ds}{l^{(1)}_n(-q_{k,n})}\right|^2
\leqslant
\sum_{n=1}^{\infty}\left(\sum\limits_{k=1}^{N-1}\left|
\frac{(-q_{k,n})\int\limits_0^tu_n(s)e^{-q_{k,n}(t-s)}ds}{l^{(1)}_n(-q_{k,n})}\right|\right)^2\leqslant
$$
$$
C_8\sum_{n=1}^{\infty}\frac{1}{\alpha_n^4}\int\limits_0^t|u_n(s)|^2ds\leqslant
C_8\sum_{n=1}^{\infty}\int\limits_0^T|u_n(s)|^2ds.
$$
Thus, uniform in $t\in[0,T]$ convergence in the norm $L_2(\Omega)$ of the
series for $\theta_t$ is proved.

\section{Conclusion}

The article is devoted to the problem of  distributed (throughout the domain) controllability of multi-dimensional wave equation with integral memory. The aim is to drive the system to rest by means of control function which is bounded by its absolute value.
The additional difficulties deal with the fact that the driving of a solution (and its first derivative, with respect to time) to null is not equivalent to controllability to rest. More accurately, not every control function, which drives the solution and its time derivative to null, leaves them in this state in the future.
In order to satisfy the last condition $u(t,x)$ should be restricted by the additional requirements (\ref{18}).

Note that when control function is applied only to the sub-domain, controllability to rest is
(generally speaking) impossible (see \cite{JOTA}).

\end{document}